\documentclass[12pt]{article}

\usepackage{amsmath,amsthm,amssymb}

\theoremstyle{plain} 
\newtheorem{thm}{Theorem}

\newtheorem{lem}[thm]{Lemma}
\newtheorem{prop}[thm]{Proposition}

\theoremstyle{definition}

\newtheorem{rem}{Remark}

\title{
$G$-reconstruction of graphs
\footnote{Published in Ars Combinatoria 54 (2000) 293-299.}
}
\author{Bhalchandra D. Thatte \\
\small \texttt{bdthatte@gmail.com} 
}
\date{}

\begin{document}
\maketitle{}
\begin{abstract}
Let $G$ be a group of permutations acting on an $n$-vertex set $V$, and $X$
and $Y$ be two simple graphs on $V$. We say that $X$ and $Y$ are 
$G$-isomorphic if $Y$ belongs to the orbit of $X$ under the 
action of $G$. One can naturally generalize the reconstruction 
problems so that when $G$ is $S_n$, the symmetric group, we have 
the usual reconstruction problems. In this paper, we study 
$G$-edge reconstructibility of graphs. 
We prove some old and new results on edge reconstruction and  
reconstruction from end vertex deleted subgraphs.  
\end{abstract}

\section{Introduction} 
\label{intro}
Unless specified, 
all the graphs in this paper are 
assumed to be undirected and without multiedges or loops, and 
to have $n$ vertices and $m$ edges.
Distance between any two vertices $u$ and $v$ is denoted by 
$d(u,v)$, and maximum degree in  a graph $X$ is denoted by $\Delta (X)$ or
simply $\Delta $ when there is no confusion. 
The automorphism group of a graph $X$ is denoted by aut$X$.

The vertex deck of a graph $X$, denoted by $VD(X)$, is the collection of 
all its unlabelled vertex deleted subgraphs, and the graph $X$ (or a 
property or a parameter of $X$) is said to be vertex reconstructible if 
$X$ (or the property or the parameter) can be uniquely obtained 
from $VD(X)$. Similarly one also defines edge deck $ED(X)$ and edge 
reconstructibility. The collection of unlabelled subgraphs of $X$ obtained by 
deleting degree-1-vertices, called end-vertex deck, 
is denoted by $VD_1(X)$, and correspondingly we
have end-vertex reconstructi\-bi\-lity.   
Vertex reconstruction conjecture (VRC) states that 
graphs with at least three vertices are vertex 
reconstructible. Edge reconstruction conjecture states that
graphs with at least four edges are edge reconstructible.
One can also pose the same conjectures in the language of hypomorphisms
between labelled graphs as follows. Two graphs $X$ and $Y$ are 
said to be vertex 
hypomorphic, denoted by $X\sim Y$, when there is a bijection $f$,
called vertex hypomorphism, from $V(X)$ to $V(Y)$ such that 
$X-u\cong Y-f(u)$ for all $u \in V(X)$. VRC then states 
that $X\sim Y$ implies $X\cong Y$, provided $n\geq 3$. 
One 
similarly defines edge hypomorphism and can pose ERC.
Reader is referred to [B] for survey of various reconstruction 
problems.

Let $V(X)=V(Y)=V$, and $G$ be a group 
of permutations acting on $V$. The action of $G$ defines the orbits of 
$X$ and $Y$, and we say that $X$ and $Y$ are $G$-isomorphic,
denoted by $X\stackrel{  G}{\cong} Y$ if $Y$ is in the orbit 
of $X$ under the action of $G$.
We can then naturally extend the above definitions to $G$-vertex (or edge)
hypomorphism, (denoted by the symbol $\stackrel{  G}{\sim}$), 
$G$-vertex (or edge) reconstructibility etc., and study the 
corresponding reconstruction problems.

Given a graph $X$ and and an edge set $E \subseteq
E(X)$, an edge set $F$ is called a 
replacing edge set of $E$ if $X-E+F\sim X$ 
(or $X-E+F\stackrel{  G}{\sim}X$)
and $E\cap F
=\emptyset $.

In this paper we demonstrate that edge or vertex reconstructibility 
of graphs can be proved under some circumstances by suitably 
choosing a group $G$ and considering the problem as $G$-ERC or $G$-VRC.
In Section~\ref{tool}, we state a generalization of 
the well known Nash-Williams' lemma. 
It is then applied to ERC in Section~\ref{edge}, and to vertex reconstruction 
of graphs from their end-vertex deleted subgraphs in Section~\ref{vertex}. 

This is an expanded version of [T3].

\section{$G$-edge reconstruction}
\label{tool}

Let $V(X)=V(Y)=V$ and $F$ be a spanning subgraph of $X$. 
For a group $G$, we denote 
$|\{g\in G|g(Y)\cap X =F\}|$ by 
$|Y\stackrel{  G}{\longrightarrow} X|_F$.
The following lemma, which is a generalization of the Nash-Williams'
lemma, is our tool in dealing with the reconstruction 
problems considered in the next two subsections.
It can be proved along the same lines as Theorem 1.1 in [T1],
and also follows from Theorem 2.1 in [ACKR].

\begin{lem}
\label{nw}
If $X$ and $Y$ are $G$-edge hypomorphic but not $G$-isomorphic then 
for every spanning subgraph $F$ of $X$, we have 
\[
|X\stackrel{  G}{\longrightarrow} X|_F
-
|Y\stackrel{  G}{\longrightarrow} X|_F
=
(-1)^{m-|E(F)|}|G\cap \rm{aut}X|
\]
\end{lem}

In the following, we demonstrate that many reconstruction problems 
can be naturally formulated as $G$-edge 
reconstruction problems, and 
Lemma~\ref{nw} can be applied.

\subsection{Edge reconstruction}
\label{edge}

The graphs considered in this section are $2$-edge connected 
bipartite graphs or separable graphs with $3$-connected pruned centre.

\subsubsection*{$2$-edge connected bipartite graphs}

\begin{prop}
\label{bipartite}
Let $X$ be a  
$2$-edge connected bipartite graph with $s$ and $t$
as the sizes of the two parts. Then $X$ is edge reconstructible 
provided $m > st/2$ or $2^{m-1}> |\rm{aut}K_{s,t}|$.  
\end{prop}

\begin{proof} 
The recognition is trivial. Also, because of $2$-edge connectivity,
the vertex partitions are uniquely recognized in the subgraphs, so
we assume the edge hypomorphism between $X$ and a possible reconstruction $Y$
to be a $G$-edge hypomorphism, where $G$ is $\rm{aut}K_{s,t}$. Now the 
claim is a simple corollary of Lemma~\ref{nw}. 
\end{proof}
That  $m > st/2$ is sufficient for edge reconstructibility,  
was proved in [VY1].

\subsubsection*{Separable graphs with $3$-edge connected pruned centre}

For a graph $X$, define the {\em block-cutpoint tree} $T(X)$, 
whose vertex set has all the cutpoints and all the maximal two
connected subgraphs ({\em $2$-blocks}) in it, and two vertices
of $T(X)$ are joined by an edge {\em iff} one of them is a
$2$-block and other is a cutpoint on the same $2$-block. The
{\em pruned graph} $P(X)$ is the maximal subgraph without
end-vertices. The center of $T(P(X))$ corresponds to 
a $2$-block or a cut vertex of $X$, and is called the {\em pruned
center} of $X$, denoted by $C(X)$.
\begin{prop}
Let $X$ be a separable  graph
with end vertices, having a $3$-connected pruned center $C(X)$. 
Suppose we colour the vertices 
of $C(X)$ blue and vertices outside $C(X)$ red.  
Let $G'$ denote the automorphism group of the coloured graph $X-E(C(X))$,
and $G$ its subgroup induced by $V(C(X))$. Then $X$ is edge 
reconstructible if $C(X)$ is $G$-edge reconstructible.
\end{prop}

\begin{proof} 
When an edge incident 
with an end vertex is deleted, (which is easily recognizable), we know the 
pruned graph and the pruned centre uniquely. This allows us to recognize 
the subset $S$ of $ED(X)$ obtained by deleting an edge of 
$C(X)$ -- given $X-e$, 
$e\in E(C(X))$ {\em iff} $T(P(X))\cong T(P(X-e))$ and $|E(C(X))|=
|E(C(X-e))|+1$. Once $S$ is recognized, we can assume that the pruned 
centres of $X$ and 
$Y$, (where $Y\sim X$), have the same vertex set, and we are given 
only the graphs in $S$. Therefore, $C(X)
\stackrel{G}{\cong}C(Y)$ is enough for an isomorphism between 
$X$ and $Y$. 
\end{proof}

Following are some of the immediate consequences:
\begin{enumerate}
\item If $2^{|E(C(X))|} > |G| $ or  $|E(C(X))| > |V(C(X))|/2 $ then 
$X$ is edge reconstructible. Note that 
we have  $2^{|E(C(X))|}$ rather than  $2^{|E(C(X))|-1}$ 
because we know $C(X)$ uniquely. Version of Lov\'asz's result was 
proved earlier in [VY].
\item If $C(X)$ is claw-free or $P_4$-free or chordal, then $X$ is 
edge reconstructible. (A graph is chordal if no induced subgraph on 
four or more vertices is a cycle. A graph is claw-free if no induced
subgraph is isomorphic to $K_{1,3}$. A graph is $P_4$-free if no induced
subgraph on four vertices is a path.)

We in fact have something stronger: 

\begin{enumerate}
\item One can observe that all connected claw-free graphs other than paths
are $G$-edge reconstructible for all groups $G$
(irrespective of their connectivity). We do not give the details 
of the proof here, but 
refer the reader to [T2], where it is proved that a collection of
connected claw-free graphs can be reconstructed from its  shuffled edge deck.  
We only comment that all the steps in that proof 
are actually based on the fact that some edge set 
in a claw-free graph has no replacement unless it is a path.  
Paths are not $G$-reconstructible for some groups, (for example,
a $2k$-vertex path, for $k\geq 2$, is not $A_{2k}$-edge reconstructible,
where $A_{2k}$ is the alternating group).
\item In case of chordal graphs, we again follow the proofs in Section 4 of 
[T2], and claim that all $2$-connected chordal 
graphs are $G$-edge reconstructible for 
all groups $G$. We also point out that, all trees except the thirteen
trees listed in [CS] are $G$-edge reconstructible, since it is proved in 
[CS] that some edge sets in all the other trees have no replacing edge sets.

\item Connected $P_{4}$-free graphs are $G$ reconstructible for all 
groups $G$ --  
the complement of any connected $P_{4}$-free graph is 
disconnected, therefore, set of all the
edges cannot be replaced.

\end{enumerate}
\end{enumerate}

\noindent {\bf Question} Can one reduce ERC for graphs with 
$3$-connected pruned centre to $3$-connected graphs?  

\subsection{End vertex deleted subgraphs}
\label{vertex}

Let $Z$ be a graph with minimum degree at least $2$, and $X$ be a 
graph obtained by adding some more vertices of degree $1$, making the 
new vertices adjacent to vertices in $Z$. 
We give simple proofs of some results on reconstruction from
end vertex deleted 
subgraphs, some of which appeared in [L].
First we consider the case in which $r_1$ vertices of $Z$ are made 
adjacent to $r_1$ new vertices of degree $1$.

\begin{prop}
\label{r1}
If $r_1 > |V(Z)|/2 $ or $2^{r_1-1}> \rm{aut}Z$
then $X$ can be reconstructed from its end-vertex 
deleted subgraphs.
\end{prop}

\begin{proof} 
Since there are no end vertices in $Z$, vertices of $Z$ are recognizable 
in every member of the deck. 
We identify all the end vertices of $X$,
call the resulting vertex $v$ and colour it blue, and rest of the 
vertices red. Call this graph $Y$. Since no two end vertices of $X$
have a common neighbour, $X$ is reconstructible from its end vertex deleted 
subgraphs, if $Y$ is reconstructible from its subgraphs $Y-av$, where $av$ is 
an edge between $v$ and a red vertex $a$. Thus we are just edge 
reconstructing a graph which is a disjoint union of an $r_1$-star and 
some isolated vertices, with the centre of the star coloured blue, with 
respect to the group aut$Z$. A direct application of 
Lemma~\ref{nw} gives the result.
\end{proof}

Now, we extend this idea to prove something stronger.
Let $Z$ be as above. Add end vertices in this graph to 
construct $X$ as follows. Let $R_i\subseteq V(Z)$, $i=1$ to $k$ be 
disjoint sets, and $|R_i|=r_i\geq 0$. For $i=1$ to $k$, we join 
each member of $R_i$ to precisely $i$ end-vertices. 
Set of remaining vertices of $Z$ is denoted by $R_0$.   
Following result is somewhat stronger than the results in [L].
\begin{prop}
\label{ri}
If $r_j>r_{j-1}$ or $2^{r_j-1}> |\rm{aut}Z|$ 
for some $j \leq k$, then $X$ is end-vertex reconstructible.
\end{prop}
\begin{proof}
Let $VD_1(Y)=VD_1(X)$. 
First we make some `recognition' claims. 
As in Proposition~\ref{r1}, vertices of $Z$ are recognized 
in every graph in the deck.
Also, the $P(Y)\cong  Z$. If $Y$ is obtained 
as above by joining $r_i'$ vertices to $i$ end-vertices each, then 
$r_i'=r_i$. This is trivial to prove unless 
$r_1=2$, $r_i=0$ for $i\geq 2$ and,  $r_1'=0$, $r_2=1$ and 
$r_i'=0$ for $i\geq 3$, in which case there are obvious counter examples.
It is also trivial to recognize, for any given graph in the deck,
the $i$ for which a neighbour of $u\in R_i$ is deleted.    
Thus we have a natural partition of the given deck into 
decks $D_i$, $i=1$ to $k$, where a member of $D_i$ results 
from deleting an end-vertex adjacent to a vertex in $R_i$.
Obviously, the multiplicity of any unlabelled graph in $D_i$
is a multiple of $i$, so we can construct reduced decks $D_i'$
by reducing all the multiplicities by a factor of $i$. 

Construct a graph $X_j$, $j\geq 1 $ from $Z$ 
by colouring vertices of $R_i$ with colour $c_i$, where 
$i\not\in \{j-1,j\}$ and making vertices of $R_j$ adjacent to 
one end-vertex each. If $X_j$ is end-vertex reconstructible, then
$X$ is end-vertex reconstructible. Therefore,
if  $r_j>r_{j-1}$ or $2^{r_j-1}> |\rm{aut}X_j|$
then $X$ is end-vertex reconstructible,  
as in Proposition~\ref{r1}. Note that $|\rm{aut}X_j| \leq |\rm{aut}Z|$. 

For some $j$, if $r_j> |V(Z)|/2$ then the graph is 
end-vertex reconstructible, which was proved in [La].
This is a corollary of our results in Propositions 2.4 and 2.5.    
\end{proof}

\begin{rem} Suppose that 
we are given the vertex deck of an arbitrary
separable graph $X$ with 
end-vertices. Given any end-vertex deleted 
graph $X-u$, it is easy 
to recognize the distance of $u$ from the nearest 
vertex of $P(X)$. Thus we know the number of end vertices at 
any distance $j$ from $P(X)$, and we can prove analogous results
on the vertex reconstruction of $X$. 
\end{rem}
  
\section*{Acknowledgements}
I did this work at the Department of Mathematics, 
Indian Institute of \text{Science}, Bangalore, India,
and at Combinatorics and Optimization, University of \text{Waterloo}, 
Canada. In Bangalore, I was supported by a post-doctoral fellow\-ship of 
the National Board for Higher Mathematics, India. At Waterloo, I
was supported by Adrian Bondy under the grant NSERC A7331.
I would like to thank Adrian Bondy for supporting my visit and for his
encouragement, and other members at C\&O, Waterloo for their help during
my visit.
 
\section*{References}  
\begin{itemize}
\item[{[ACKR]}] N. Alon, Y. Caro, I. Krasikov and Y. Roditty, Combinatorial
reconstruction problems, {\em J. Combin. Theory Ser. B} {\bf 47} no. 2 (1989)
153-161. 
\item[{[B]}]    J. A. Bondy, A graph reconstructor's manual, in "Surveys in
Combinatorics, 1991", (Guildford, 1991), London Math. Soc. Lecture Note
Ser. 166, Cambridge Univ. Press, Cambridge (1991) 221-252.
\item[{[CS]}]  C. R. J. Clapham and J. Sheehan, The thirteen two-free trees,
preprint, Dept of Mathematical Sciences, Univ. of Aberdeen (1991).
\item[{[La]}]  J. Lauri, End vertex deleted subgraphs,
{\em Ars Combinatoria}, {\bf 36} (1993), 171-184. 
\item[{[Lo]}] L. Lov\'asz, A note on the line reconstruction problem, 
{\em J. Combin.
Theory Ser. B} {\bf 13} (1972) 309-310. 
\item[{[Mn]}] V. B. Mnukhin, The $k$-orbit reconstruction and the orbit
algebra, {\em Acta Applicandae Mathematicae} {\bf 29} (1992) 83-117.
\item[{[M\"u]}] V. M\"uller, The edge reconstruction hypothesis is true for
graphs with more than $n.\log_2 n$ edges, {\em J. Combin. Theory Ser. B}
{\bf 22} no. 3 (1977) 281-283. 
\item[{[T1]}] B.D. Thatte, On the Nash-Williams'lemma in graph
reconstruction theory, {\em J.
Combin. Theory Ser. B} {\bf 58},2,(1993)280-290.
\item[{[T2]}] B.D. Thatte, Some results and approaches for
reconstruction problems, The first Malta Conference on Graphs and
Combinatorics, May-June 1990, {\em Discrete Mathematics}, {\bf 124}(1994),
193-216.
\item[{[T3]}]  B. D. Thatte, More appli\-cations of a gene\-ral Nash-\-Williams'
lemma, manu\-script, (1992).
\item[{[VY1]}] A. Vince and Y. Yong Zhi, Complement edge
reconstruction, Preprint, (1992).
\item[{[VY1]}] A. Vince and Y. Yong Zhi, Edge
reconstruction of graphs with sufficiently large center, Preprint, (1992).
\end{itemize}
\end{document}